\documentstyle[12pt,amssymb,amscd]{amsart}
\chardef\bslash=`\\ 

\makeatletter
\def\verbatim{\interlinepenalty\@M \@verbatim
  \leftskip\@totalleftmargin\advance\leftskip2pc
  \frenchspacing\@vobeyspaces \@xverbatim}
\renewcommand{\thesubsection}{\thesection(\@roman\c@subsection)}
\makeatother
\newtheorem{Theorem}[equation]{Theorem}
\newtheorem{Corollary}[equation]{Corollary}
\newtheorem{Lemma}[equation]{Lemma}
\newtheorem{Proposition}[equation]{Proposition}

\theoremstyle{definition}

\theoremstyle{remark}
\newtheorem{Remark}[equation]{Remark}

\numberwithin{equation}{section}

\newcommand{\thmref}[1]{Theorem~\ref{#1}}
\newcommand{\secref}[1]{\S\ref{#1}}
\newcommand{\lemref}[1]{Lemma~\ref{#1}}
\newcommand{\subsecref}[1]{\S\ref{#1}}

\newcommand{\defeq}{\overset{\operatorname{\scriptstyle def.}}{=}}
\newcommand{\interior}{\rule{6pt}{.4pt}\rule{.4pt}{6pt}} 
\newcommand{\C}{{\Bbb C}}

\newcommand{\R}{{\Bbb R}}

\newcommand{\CP}{\operatorname{\C P}}

\newcommand{\SU}{\operatorname{\rm SU}}
\newcommand{\GL}{\operatorname{GL}}
\newcommand{\U}{\operatorname{\rm U}}

\newcommand{\algsl}{\operatorname{\frak{sl}}} 

\newcommand{\Hom}{\operatorname{Hom}}

\newcommand{\Ker}{\operatorname{Ker}}

\newcommand{\Ima}{\operatorname{Im}}

\newcommand{\rank}{\operatorname{rank}}

\newcommand{\tr}{\operatorname{tr}}

\newcommand{\Hilb}[2]{{#1}^{\lbrack{#2}\rbrack}} 
\newcommand{\Hilbn}[1]{\Hilb{#1}{n}}
\newcommand{\HilbX}[1]{\Hilb{X}{#1}}
\newcommand{\shfO}{\cal O} 
\newcommand{\Hlf}{H^{\it lf}} 
\newcommand{\Cl}{\operatorname{Cl}} 
\newcommand{\idl}{\cal J} 
\newcommand{\Supp}{\operatorname{Supp}} 
\newcommand{\codim}{\operatorname{codim}}
\newcommand{\bM}{{\bold M}} 
\newcommand{\M}{{\frak M}} 

\begin{document}
\title[Heisenberg Algebra and Hilbert Schemes]
{Heisenberg Algebra and Hilbert Schemes of Points on Projective Surfaces}
\author{Hiraku Nakajima}
\address{Department of Mathematical Sciences, University of Tokyo\\
Komaba 3-8-1, Meguro-ku, Tokyo 153, Japan}
\email{nakajima@@math.tohoku.ac.jp}
\thanks{Author supported in part by Grant-in-Aid for Scientific Research
(No.\ 05740041), Ministry of Education, Science and Culture, Japan
and also by the Inamori Foundation.
}
\maketitle

%
%
\section{Introduction}
The purpose of this paper is to throw a bridge between two seemingly
unrelated subjects. One is the Hilbert scheme of points on projective
surfaces, which has been intensively studied by various people (see
e.g., \cite{Iar,ES,Got,Go-book}). The other is the infinite
dimensional Heisenberg algebra which is closely related to affine Lie
algebras (see e.g., \cite{Kac}).

We shall construct a representation of the Heisenberg algebra on the
homology group of the Hilbert scheme.  In other words, the homology
group will become a Fock space.  The basic idea is to introduce
certain ``correspondences'' in the product of the Hilbert scheme.
Then they define operators on the homology group by a well-known
procedure.  They give generators of the Heisenberg algebra, and the
only thing we must check is that they satisfy the defining relation.
Here we remark that the components of the Hilbert scheme are
parameterized by numbers of points and our representation will be
constructed on the direct sum of homology groups of all components.
Our correspondences live in the product of the different components.
Thus it is quite essential to study all components together.

Our construction has the same spirit with
author's construction \cite{Na-quiver,Na-gauge}
of representations of affine Lie algebras
on homology groups of moduli spaces of ``instantons''\footnote{The
reason why we put the quotation mark will be
  explained in Remark~\ref{sheaf}.}
on ALE spaces which are minimal resolution of simple singularities.
Certain correspondences, called Hecke correspondences, were used to
define operators. These twist vector bundles along
curves (irreducible components of the exceptional set), while ours
twist around points. In fact, the Hilbert scheme of points can be
considered as the moduli space of rank $1$ vector bundles, or more
precisely torsion free sheaves.
Our construction should be considered as a first step to extend
\cite{Na-quiver,Na-gauge} to more general $4$-manifolds.
The same program was also proposed by Ginzburg, Kapranov and Vasserot
\cite{GKV}.

Another motivation of our study is the conjecture about the generating
function of the Euler number of the moduli spaces of instantons,
which was recently proposed by Vafa and Witten~\cite{VW}.
They conjectured that it is a modular form for $4$-manifolds under
certain conditions.
This conjecture was checked for various $4$-manifolds using
various mathematicians' results.
Among them, the most relevant to us is the case of K3 surfaces.
G\"ottsche and Huybrechts \cite{GotHu} proved that the Betti numbers
of moduli spaces of stable rank two sheaves are the same as those for
Hilbert schemes. G\"ottsche \cite{Got} computed the Betti numbers
of Hilbert schemes for general projective surfaces $X$.
(The Hilbert schemes for $\CP^2$ were studied earlier by Ellingsrud and
Str\o mme's~\cite{ES}.)
If $\Hilbn{X}$ is the Hilbert scheme parameterizing $n$-points in $X$,
the generating function of the Poincar\'e polynomials is given by
\begin{equation}
  \sum_{n=0}^\infty q^n P_t(\Hilbn{X}) = \prod_{m=1}^\infty
  \frac{(1 + t^{2m-1}q^m)^{b_1(X)}(1 + t^{2m+1}q^m)^{b_3(X)}}
       {(1 - t^{2m-2}q^m)^{b_0(X)}(1 - t^{2m}q^m)^{b_2(X)}
                 (1 - t^{2m+2}q^m)^{b_4(X)}}\, ,
\label{Poincare}\end{equation}
where $b_i(X)$ is the Betti number of $X$. Letting $t = -1$, we find
the generating function of the Euler numbers is essentially the
Dedekind eta function. In fact, the relation with the above formula
and the Fock space was already pointed out in \cite{VW}. Our result
should be considered as a geometric realization of their indication.

The paper is organized as follows.
In \secref{sec:pre} we give preliminaries. We recall the definition of
the convolution product in \subsecref{subsec:conv}
with some modifications and describe some
properties of the Hilbert scheme $\HilbX{n}$ and the infinite
Heisenberg algebra and its representations
\S\S\ref{subsec:Hilb},\ref{subsec:Heisen}.
The definition of correspondences and the statement of the main result
are given in \secref{sec:main}.
The proof will be given in \secref{sec:proof}.
In the appendix, we study the particular case $X = \C^2$ in more
detail.
We give a description
of $\Hilbn{(\C^2)}$ as a hyper-K\"ahler quotient of
finite dimensional vector space by a unitary group action.
It is very similar to the definition of quiver varieties
\cite{Na-quiver}. The only difference is that we have an edge joining
a vertex with itself.
Using this description as a hyper-K\"ahler quotient, we compute the
homology group of $\Hilbn{(\C^2)}$.
We recover the formula \eqref{Poincare} for $X = \C^2$.
The difference between our approach and
Ellingsrud-Str\o mme's~\cite{ES} is only the description.
Both use the torus action and study the fixed point set.
But our presentation has a similarity in \cite{Na-homology}.
The appendix is independent of the other parts of this paper, but
those similarities with author's previous works explains motivation of
this paper in part.

While the author was preparing this manuscript, he learned that the
similar result was announced by Grojnowski \cite{Gr}. He introduced
exactly the same correspondence as ours.

\subsection*{Acknowledgement}
The author would like to thank C.~Vafa and E.~Witten, since
it is clear that this work was not done unless he discussed
with them.
It is also a pleasure to acknowledge discussions with V.~Ginzburg and
M.~Kapranov.
His thanks go also to R.~Hotta, T.~Uzawa, K.~Hasegawa and G.~Kuroki
who answered many questions on the representation theory.
\section{Preliminaries}
\label{sec:pre}
\subsection{Convolution Algebras}
\label{subsec:conv}
We need a slight modification of the definition of the convolution
product in the homology groups given by Ginzburg \cite{Gi} (see also
\cite{Gi-book,Na-quiver}).

For a locally compact topological space $X$, let $\Hlf_*(X)$
denote the homology group of possibly infinite singular chains with locally
finite support (the Borel-Moore homology) with {\it rational\/}
coefficients.
The usual homology group of finite singular chains will be denoted by
$H_*(X)$.
If $\overline X = X\cup \{\infty\}$ is the one point compactification
of $X$, we have $\Hlf_*(X)$ is isomorphic to the relative homology
group $H_*(\overline X, \{\infty\})$.
If $X$ is an $n$-dimensional oriented manifold,
we have the Poincar\'e duality isomorphism
\begin{equation}
        \Hlf_i(X) \cong H^{n-i}(X),\quad
        H_i(X) \cong H^{n-i}_c(X),
\label{eq:PD}\end{equation}
where $H^*$ and $H^*_c$ denote the ordinary cohomology group and the
cohomology group with compact support respectively.

%
%
%

Let $M^1$, $M^2$, $M^3$ be oriented manifolds of dimensions $d_1$,
$d_2$, $d_3$ respectively, and
$p_{ij}\colon M^1\times M^2\times M^3 \to M^i\times M^j$ be the
natural projection.
We define a convolution product
\begin{equation*}
  \ast\colon \left(\Hlf_{i_1}(M^1)\otimes H_{i_2}(M^2)\right) \otimes
  \left(\Hlf_{d_2 - i_2}(M^2)\otimes H_{i_3}(M^3)\right) \to
  \Hlf_{i_1}(M^1)\otimes H_{i_3}(M^3)
\end{equation*}
by
\begin{equation*}
  (c_1\otimes c_2) \ast (c'_2\otimes c_3) \defeq
  c_2\cap c'_2\; c_1\otimes c_3,
\end{equation*}
where $c_2\cap c'_2\in\Bbb Z$ is the natural pairing between
$H_{i_2}(M^2)$ and $\Hlf_{d_2 - i_2}(M^2)\cong H^{i_2}(M^2)$.

Suppose $Z$ is a submanifold in $M^1\times M^2$ such that
\begin{equation}
  \text{the projection $Z \to M^1$ is proper.}
  \label{ass-proper}
\end{equation}
Then the fundamental class $[Z]$ defines an element in
\begin{equation*}
  [Z] \in
  H_{\dim_\R Z}(\overline M_1\times M_2, \{\infty\}\times M_2)
  = \bigoplus_{i + j = \dim_\R Z} \Hlf_i(M^1)\otimes H_j(M^2),
\end{equation*}
where $\overline M_1 = M_1\cup\{\infty\}$ is the one point
compactification of $M_1$ and we have used the K\"unneth formula.
More generally, if $[Z]$ is a cycle whose support $Z$ satisfies
\eqref{ass-proper}, the same construction works.
Using \eqref{eq:PD}, we get an operator, which is denoted also by $[Z]$,
\begin{equation*}
  [Z]\colon \Hlf_j(M^2)\to \Hlf_{j + \dim_\R Z - d_2}(M^1).
\end{equation*}
\subsection{Hilbert Schemes of Points on Surfaces}
\label{subsec:Hilb}
Let $X$ be a nonsingular quasi-projective surface defined
over the complex number $\C$. Let $\Hilbn{X}$ be the component of the
Hilbert scheme of $X$ parameterizing the ideals of $\shfO_X$ of
colength $n$. It is smooth and irreducible \cite{Fogarty}.
Let $S^n X$ denotes the $n$-th symmetric product of $X$.
It parameterizes formal linear combinations $\sum n_i [x_i]$ of points
$x_i$ in $X$ with coefficients $n_i\in\Bbb Z_{> 0}$ with $\sum n_i = n$.
There is a canonical morphism
\begin{equation*}
  \pi\colon \Hilbn{X}\to S^n X; \quad
  \pi(\cal J) \defeq \sum_{x\in X}
           \operatorname{length}(\shfO_X/\cal J)_x [x].
\end{equation*}
It is known that $\pi$ is a resolution of singularities.

The symmetric power $S^n X$ has a natural stratification into locally
closed subvarieties as follows. Let $\nu$ be a partition of $n$, i.e.,
a sequence $n_1, n_2, \dots, n_r$ such that
\begin{equation*}
  n_1 \ge n_2 \ge \cdots \ge n_r, \qquad \sum n_i = n.
\end{equation*}
Then $S^n_\nu X$ is defined by
\begin{equation*}
  S^n_\nu X \defeq \{ \sum_i n_i [x_i] \in S^n X\mid
     x_i \ne x_j \quad\text{for $i\ne j$}\, \}.
\end{equation*}

It is known that $\pi$ is semi-small with respect to the
stratification $S^n X = \bigcup S^n_\nu X$ \cite{Iar}, that is
\begin{enumerate}
\item for each $\nu$, the restriction
$\pi\colon \pi^{-1}(S^n_\nu X)\to S^n_\nu X$ is a locally trivial
fibration,
\item $\codim S^n_\nu X = 2\dim \pi^{-1}(x)$ for $x\in S^n_\nu X$.
\end{enumerate}
Moreover, it is also known that $\pi^{-1}(x)$ is irreducible.
\subsection{The Infinite Dimensional Heisenberg Algebra}
\label{subsec:Heisen}
We briefly recall the definition of the infinite dimensional
Heisenberg algebra and its representations. See \cite[\S9.13]{Kac} for
detail.

The infinite dimensional Heisenberg algebra $\frak s$ is generated by
$p_i$, $q_i$ ($i=1,2,\dots$) and $c$ with the following relations:
\begin{gather}
  [p_i, p_j] = 0, \quad [q_i, q_j] = 0 \label{eq:rel1}\\
  [p_i, q_j] = \delta_{ij} c. \label{eq:rel2}
\end{gather}

For every $a\in\C^*$, the Lie algebra $\frak s$ has an irreducible
representation on the space $R = \C[x_1, x_2,\dots]$ of
polynomials in infinitely many indeterminates $x_i$ defined by
\begin{equation*}
  p_i \mapsto a \frac{\partial}{\partial x_i}, \quad
  q_i \mapsto x_i, \quad
  c \mapsto a \operatorname{Id}.
\end{equation*}
This representation has a highest weight vector $1$, and $R$ is
spanned by elements
\begin{equation*}
  x_1^{j_1} x_2^{j_2} \cdots x_n^{j_n} =
  q_1^{j_1} q_2^{j_2} \cdots q_n^{j_n} 1.
\end{equation*}

We extend $\frak s$ by a derivation $d_0$ defined by
\begin{equation*}
  [d_0, q_j] = j q_j, \quad [d_0, p_j] = - j p_j.
\end{equation*}
The above representation $R$ extends by
\begin{equation*}
  d_0 \mapsto \sum_j jx_j\frac{\partial}{\partial x_j}.
\end{equation*}
Then it is easy to see
\begin{equation}
  \tr_R q^{d_0} = \prod_{j=1}^\infty \frac 1{(1 - q^j)}.
  \label{eq:char}
\end{equation}

The representation $R$ carries a unique bilinear form $B$ such that
$B(1, 1) = 1$ and $p_i$ is the adjoint of $q_i$, provided $a\in \R$.
In fact, distinct monomials are orthogonal and we have
\begin{equation*}
  B(x_1^{j_1}\dots x_n^{j_n}, x_1^{j_1}\dots x_n^{j_n})
  = a^{\sum j_k} \prod j_k !
\end{equation*}

We also need the infinite dimensional Clifford algebra $\Cl$
(see e.g., \cite{Fr}). It
is generated by $\psi_i$, $\psi^*_i$ ($i=1,2,\dots$) and $c$ with the
relations
\begin{gather}
  \psi_i\psi_j + \psi_j \psi_i = 0, \quad
  \psi_i^*\psi_j^* + \psi_j^* \psi_i^* = 0\label{eq:rel3}\\
  \psi_i\psi_j^* + \psi_j^*\psi_i = \delta_{ij}c.\label{eq:rel4}
\end{gather}
This algebra has a representation on the exterior algebra
$F = \bigwedge^* V$ of an infinite dimensional vector space
$V = \C dx^1\oplus \C dx^2\oplus\cdots$ defined by
\begin{equation*}
  \psi_i \mapsto dx^i \wedge\, , \quad
  \psi_i^* \mapsto \frac{\partial}{\partial x_i}\interior\, ,\quad
  c \mapsto \operatorname{Id},
\end{equation*}
where $\interior$ denotes the interior product.
This has the highest weight vector $1$ and spanned by
\begin{equation*}
  dx_{i_1}\wedge\cdots\wedge dx_{i_n} =
  \psi_{i_1}\cdots \psi_{i_n} 1, \qquad (i_1 > i_2 > \cdots > i_n).
\end{equation*}
We extend $\Cl$ by $d$ defined by
\begin{equation*}
  [d, \psi_i] = i\psi_i, \quad [d, \psi_i^*] = -i\psi_i.
\end{equation*}
It acts on $F$ by
\begin{equation*}
  d(dx_{i_1}\wedge\cdots\wedge dx_{i_n})
  = \left(\sum i_k\right) dx_{i_1}\wedge\cdots\wedge dx_{i_n}.
\end{equation*}
The character is given by
\begin{equation*}
  \tr_F q^d = \prod_{j=1}^\infty (1 + q^j).
\end{equation*}
\section{Main Construction}
\label{sec:main}
\subsection{Definitions of Generators}
Let $X$ as in \subsecref{subsec:Hilb}.
Take a basis of $\Hlf_*(X)$ and
assume that each element is represented by a (real) closed submanifold
$C^a$. ($a$ runs over $1, 2, \dots, \dim \Hlf_*(X)$.)
Take a dual basis for $H_*(X)\cong H^{4 - *}_c(X)$, and assume
that each element is represented by a submanifold $D^a$ which is
compact.
(Those assumptions are only for the brevity. The modification to
the case of cycles is clear.)
For each $a = 1,2,\dots,\dim\Hlf_*(X)$, $n = 1,2,\dots$ and
$i=1,2,\dots$, we introduce cycles of products of the Hilbert
schemes by
\begin{equation*}
  \begin{split}
    &E_i^a(n) \defeq \{\, (\cal J_1,\cal J_2)\in\HilbX{n-i}\times\HilbX{n}
    \mid \\
    & \qquad\qquad\qquad\qquad\text{$\cal J_1\supset \cal J_2$ and
      $\Supp(\cal J_1/\cal J_2) = \{ p\}$ for some $p\in D^a$}
      \}, \\
    &F_i^a(n) \defeq \{\, (\cal J_1,\cal J_2)\in\HilbX{n+i}\times\HilbX{n}
    \mid \\
    & \qquad\qquad\qquad\qquad\text{$\cal J_1\subset \cal J_2$ and
      $\Supp (J_2/\cal J_1) = \{ p \}$ for some $p\in C^a$}\}.
  \end{split}
\end{equation*}
The dimensions are given by
\begin{align*}
  & \dim_\R E_i^a(n) = 4(n-i) + 2(i-1) + \dim_\R D^a,\\
  & \dim_\R F_i^a(n) = 4n + 2(i-1) + \dim_\R C^a.
\end{align*}
This follows from the fact $\Hilbn{X}\to S^n X$ is semi-small (see
\subsecref{subsec:Hilb}).
Since the projections
$E_i^a(n)\to \HilbX{n-i}$ and $F_i^a(n)\to\HilbX{n+i}$ are proper,
we have classes
\begin{align*}
  & [E_i^a(n)]\in \bigoplus_{k,l}
    \Hlf_k(\HilbX{n-i})\otimes H_l(\HilbX{n}), \\
  & [F_i^a(n)]\in \bigoplus_{k,l}
    \Hlf_k(\HilbX{n+i})\otimes H_l(\HilbX{n}).
\end{align*}

Our main result is the following:
\begin{Theorem}
  The following relations hold in
  $\bigoplus_{k,l,m,n} \Hlf_k(\HilbX{m})\otimes H_l(\HilbX{n})$.
  \begin{gather}
    [E_i^a(n-j)]\ast [E_j^b(n)] =
    (-1)^{\dim D^a\dim D^b}[E_j^b(n-i)] \ast [E_i^a(n)]\label{eq:EE}\\
    [F_i^a(n+j)]\ast [F_j^b(n)] =
    (-1)^{\dim C^a\dim C^b}[F_j^b(n-i)] \ast [F_i^a(n)]\label{eq:FF}\\
    [E_i^a(n+j)]\ast [F_j^b(n)] =
    (-1)^{\dim D^a\dim C^b}[F_j^b(n-i)] \ast [E_i^a(n)]
    + \delta_{ab}\delta_{ij}c_i [\Delta(n)]\label{eq:EF},
  \end{gather}
  where $\Delta(n)$ is the diagonal of $\Hilbn{X}$, and $c_i$ is a
  nonzero integer depending only on $i$ \rom(independent of $X$\rom).
\label{th:main}\end{Theorem}
In particular, for each fixed $a$, the map
\begin{align*}
  & p_i \mapsto \sum_n [E_i^a(n)], \quad q_i \mapsto \sum_n [F_i^a(n)]
  \qquad\text{when $\dim C^a$ is even} \\
  & \psi_i^*\mapsto \sum_n [E_i^a(n)], \quad
  \psi_i\mapsto \sum_n [F_i^a(n)]
  \qquad\text{when $\dim C^a$ is odd}
\end{align*}
defines a homomorphism from the Heisenberg algebra and the Clifford
algebra respectively.

Considering $[E_i^a(n)]$, $[F_i^a(n)]$ as operators on
$\bigoplus_{k,n} \Hlf_k(\HilbX{n})$, we have a representation of the
product of Heisenberg algebras and Clifford algebras. Comparing
G\"ottsche's Betti number formula and the character formula, we get
the following:
\begin{Theorem}
  The direct sum $\bigoplus_{k,n} \Hlf_k(\HilbX{n})$ of homology groups
  of $\HilbX{n}$ is the highest weight module where the highest weight
  vector $v_0$ is the generator of $\Hlf_0(\HilbX{0})\cong\Bbb Q$.
\end{Theorem}

\begin{Remark}
  The author does not know the precise values of $c_i$'s. It is easy
  to get $c_1 = 1$, $c_2 = -2$, but general $c_i$ become difficult to
  calculate.
\end{Remark}
\section{Proof of {\protect\thmref{th:main}}}
\label{sec:proof}
\subsection{Proof of Relations (I)}
Consider the product $\HilbX{n-i-j}\times\HilbX{n-j}\times\HilbX{n}$
and let $p_{12}$, etc. be as in \subsecref{subsec:conv}. The
intersection $p_{12}^{-1}(E_i^a(n-j))\cap p_{23}^{-1}(E_j^b(n))$
consists of triples $(\idl_1, \idl_2, \idl_3)$ such that
\begin{align}
  & \idl_1\supset\idl_2\supset\idl_3 \label{eq:incl}\\
  & \text{$\Supp(\idl_1/\idl_2) = \{ p\}$, \;
    $\Supp(\idl_2/\idl_3) = \{ q\}$ for some $p\in D^a$, $q\in D^b$.}
 \label{eq:supp}
\end{align}
Replacing $D^b$ by $\widetilde D^b$ in the same homology class, we may
assume
$\dim D^a\cap \widetilde D^b = \dim D^a + \dim \widetilde D^b - 4$.
(If the right hand side is negative, the set is empty.)
Let $U$ be the open set in the intersection consisting points
with $p\ne q$ in \eqref{eq:supp}.
Outside the singular points of $p_{12}^{-1}(E_i^a(n-j))$,
$p_{23}^{-1}(E_j^b(n))$, the intersection is transverse along $U$.
The complement
$p_{12}^{-1}(E_i^a(n-j))\cap p_{23}^{-1}(E_j^b(n))\setminus U$
consists of $(\idl_1, \idl_2, \idl_3)$ with \eqref{eq:incl} and
\begin{equation*}
  \text{$\Supp(\idl_1/\idl_2) = \Supp(\idl_2/\idl_3) = \{p\}$ for some
    $p\in D^a\cap \widetilde D^b$.}
\end{equation*}
Its dimension is at most
\begin{equation*}
  4n - 2i - 2j - 4 + \dim D^a + \dim D^b - 4,
\end{equation*}
which is strictly smaller than the dimension of the intersection
\begin{equation*}
  4n - 2i - 2j - 4 + \dim D^a + \dim D^b.
\end{equation*}

Now consider the product
$\HilbX{n-i-j}\times\HilbX{n-i}\times\HilbX{n}$.
The
intersection $p_{12}^{-1}(E_j^b(n-i))\cap p_{23}^{-1}(E_i^a(n))$
consists of triples $(\idl_1, \idl'_2, \idl_3)$ such that
\begin{align}
  & \idl_1\supset\idl'_2\supset\idl_3 \label{eq:incl2}\\
  & \text{$\Supp(\idl_1/\idl'_2) = \{q\}$ \;
    $\Supp(\idl'_2/\idl_3) = \{p\}$ for some $q\in D^b$, $p\in D^a$.}
  \label{eq:supp2}
\end{align}
Let $U'$ be the open set in the intersection consisting points
with $p\ne q$ in \eqref{eq:supp2}.
The intersection is again transverse along $U$ outside singular sets.
The complement
$p_{12}^{-1}(E_j^b(n-i))\cap p_{23}^{-1}(E_i^a(n))\setminus U'$
has dimension is at most
\begin{equation*}
  4n - 2i - 2j - 4 + \dim D^a + \dim D^b - 4,
\end{equation*}
which is also strictly smaller than the dimension of the intersection.

There exists a homeomorphism between $U$ and $U'$ given by
\begin{equation*}
  U\ni (\idl_1, \idl_2, \idl_3) \mapsto (\idl_1, \idl'_2, \idl_3) \in U',
\end{equation*}
where $\idl'_2$ is a sheaf such that
\begin{equation*}
  \idl_1 / \idl'_2 = \idl_2 / \idl_3, \quad
  \idl'_2 / \idl_3 = \idl_1 / \idl_2.
\end{equation*}
Such $\idl'_2$ exists since supports of $\idl_1 / \idl_2$ and
$\idl_2 / \idl_3$ are different points $p$ and $q$.

Taking account of orientations and the estimate of the dimension of
the complements, we get the relation~\eqref{eq:EE}. The proof of
\eqref{eq:FF} is exactly the same.
\subsection{Proof of Relations (II)}
The proof of \eqref{eq:EF} is almost similar to the above.
Consider the product $\HilbX{n-i+j}\times\HilbX{n+j}\times\HilbX{n}$.
The intersection $p_{12}^{-1}(E_i^a(n+j))\cap p_{23}^{-1}(F_j^b(n))$
consists of triples $(\idl_1, \idl_2, \idl_3)$ such that
\begin{align}
  & \idl_1\supset\idl_2\subset\idl_3 \label{eq:incl3}\\
  & \text{$\Supp(\idl_1/\idl_2) = \{p\}$, \;
    $\Supp(\idl_3/\idl_2) = \{q\}$ for some $p\in D^a$, $q\in C^b$.}
 \label{eq:supp3}
\end{align}
Replacing $C^b$ by $\widetilde C^b$ in the same homology class, we may
assume
$\dim D^a\cap \widetilde C^b = \dim D^a + \dim \widetilde C^b - 4$.
(If the right hand side is negative, the set is empty.)
Since $\{D^a\}$ and $\{C^b\}$ are dual bases each other, the equality
holds if and only if $a = b$.
Let $U$ be the open set in the intersection consisting points
with $p\ne q$ in \eqref{eq:supp3}.

Next consider the product $\HilbX{n-i+j}\times\HilbX{n-i}\times\HilbX{n}$.
The intersection $p_{12}^{-1}(F_j^b(n-i))\cap p_{23}^{-1}(E_i^a(n))$
consists of triples $(\idl_1, \idl_2, \idl_3)$ such that
\begin{align}
  & \idl_1\subset\idl'_2\supset\idl_3 \label{eq:incl4}\\
  & \text{$\Supp(\idl'_2/\idl_1) = \{q\}$, \;
    $\Supp(\idl'_2/\idl_3) = \{p\}$ for some $q\in C^b$, $p\in D^a$.}
 \label{eq:supp4}
\end{align}
Let $U'$ be the open set in the intersection consisting points
with $p\ne q$ in \eqref{eq:supp4}.

There exists a homeomorphism between $U$ and $U'$ given by
\begin{equation*}
  U\ni (\idl_1, \idl_2, \idl_3) \mapsto (\idl_1, \idl'_2, \idl_3) \in U',
\end{equation*}
where $\idl'_2$ is
\begin{equation*}
  (\idl_1 \oplus \idl_3)/ \{ (f, f)\mid f\in \idl_2\} .
\end{equation*}
The inverse map is given by $\idl_2 = \idl_1 \cap \idl_3$.

Let $U^c$, $U^{\prime c}$ be the complement of $U$ and $U'$ respectively.
If $(\idl_1, \idl_3)$ is in the image $p_{13}(U^c)$ or
$p_{13}(U^{\prime c})$,
then $\idl_1$ and $\idl_3$ are isomorphic outside a point
$p\in D^a\cap \widetilde C^b$.
In particular, it is easy to check
\begin{equation*}
  \dim p_{13}(U^c), \, \dim p_{13}(U^{\prime c}) \le
  4n - 2i + 2j - 4 + \dim D^a + \dim C^b,
\end{equation*}
where the right hand side is the expected dimension of the intersection.
The equality holds only if $i = j$ and $\dim D^a + \dim C^b = 4$.
Moreover, $\{D^a\}$ and $\{C^b\}$ are dual bases,
when $\dim D^a + \dim C^b = 4$,
the intersection
$D^a\cap \widetilde C^b$ is empty unless $a = b$.
Thus we have checked \eqref{eq:EF} when $i\ne j$ or $a\ne b$.

Now assume $i = j$ and $a = b$.
Then $p_{13}(U^{\prime c})$ has smaller dimension and
$p_{13}(U^c)$ is union of the diagonal $\Delta(n)$ and smaller
dimensional sets.
Hence the left hand side of \eqref{eq:EF} is a multiple of $[\Delta(n)]$.
In order to calculate the multiple, we may restrict the intersection
on the open set where
\begin{enumerate}
\item $\idl_1$ and $\idl_3$ are contained in the open
stratum $\pi^{-1}(S^n_{1,1,\dots,1})$,
\item $\Supp\shfO/\idl_1$, $\Supp\shfO/\idl_3$ do not intersect with
  $D^a\cap \widetilde C^a$.
\end{enumerate}
Then it is clear that the multiple is a constant independent of $n$ and $X$,
which we denoted by $c_i$.
The only thing left is to show $c_i \ne 0$.
We may assume $X = \C^2$ and $n = i$.
We consider the quotient of $\Hilb{(\C^2)}{i}$ devided by the action of
$\C^2$ which comes from the parallel translation.
Thus $c_i$ is equal to the self-intersection number of
$[\pi^{-1}(i[0])]$ in $\Hilb{(\C^2)}{i}/\C^2$.
Since $\pi\colon \Hilb{(\C^2)}{i}\to S^i\C^2$ has irreducible fibers,
$[\pi^{-1}(i[0])]$ is the generator of $H_{2i - 2}(\Hilb{(\C^2)}{i}/\C^2)$.
Now our assertion follows from a general result which holds for any
semi-small morphism \cite[7.7.15]{Gi-book}.
That is the non-degeneracy of the intersection form on the top degree
of the fiber.
\appendix\section{Hilbert Schemes of Points on the Plane}
\label{sec:appendix}
In this appendix, we describe the Hilbert scheme $\Hilbn{(\C^2)}$, or
more generally framed moduli spaces of torsion free sheaves on $\CP^2$,
as a hyper-K\"ahler quotient of a finite dimensional vector space with
respect to a unitary group action, and then compute its homology.
\subsection{The framed moduli space of torsion free sheaves on
${\protect\CP^2}$}
Let $V$ and $W$ be vector spaces over the complex field whose
dimensions are $n$ and $r$. Let
\begin{equation*}
  \bM \defeq \{ (B_1, B_2, i, j)\mid B_1, B_2\in\Hom(V, V),
  i\in \Hom(W, V), j\in\Hom(W,V) \}.
\end{equation*}
Consider the following complex ADHM equation:
\begin{equation*}
  \mu_\C(B_1,B_2,i,j) \defeq [B_1, B_2] + ij = 0.
\end{equation*}

Let $[z_0 : z_1 : z_2]$ be the homogeneous coordinates of $\CP^2$.
We consider the following homomorphisms of sheaves over $\CP^2$
\begin{equation*}
  V\otimes\shfO(-1) @>\sigma>> (V\oplus V \oplus W)\otimes\shfO
 @>\tau>> V\otimes\shfO(1),
\end{equation*}
where
\begin{equation*}
  \sigma = \begin{pmatrix} B_1\otimes z_0 - 1_V\otimes z_1 \\
  B_2\otimes z_0 - 1_V\otimes z_2 \\ j\otimes z_0\end{pmatrix},
  \quad \tau = \begin{pmatrix}
  - (B_2\otimes z_0 - 1_V\otimes z_2) &
  B_1\otimes z_0 - 1_V\otimes z_1 & i\otimes z_0\end{pmatrix}.
\end{equation*}
Here $1_V$ denotes the identity map of $V$.
By the complex ADHM equation, this is a complex, that is
$\tau\sigma = 0$.

If $\sigma$ is injective and $\tau$ is surjective as sheaf
homomorphisms, then $E \defeq \Ker\tau / \Ima\tau$ is a torsion free sheaf on
$\CP^2$ with
\begin{equation*}
  \rank E = r, \quad c_1(E) = 0, \quad c_2(E) = n.
\end{equation*}
Remark that the injectivity of $\sigma$ is weaker than the requirement
that $\sigma$ induces injective homomorphisms between stalks over any
points.
This stronger condition makes the resulting $E$ {\it locally free\/}.

Over the line $l = \{z_0 = 0\}$, the sheaf $E$ is naturally identified
with $W\otimes\shfO_l$.

Conversely suppose we are given a torsion free sheaf $E$ over $\CP^2$
which has a trivialization $E|_l \cong W\otimes\shfO_l$.
We then define a vector space $V$ by
\begin{equation*}
  V \defeq H^1(\CP^2, E(-2)).
\end{equation*}
Using the vanishing of $H^0(l, E(-1))$ and $H^1(l, E(-1))$, we obtain
the isomorphism $H^1(\CP^2, E(-2))\cong H^1(\CP^2, E(-1))$. Thus the
multiplications of $z_1$, $z_2$ define endomorphisms $B_1$, $B_2$ of
$V$. We define $i$ and $j$ as natural homomorphisms
\begin{equation*}\begin{split}
  i &\colon W\cong H^0(l, E) \to V\cong H^1(\CP^2, E(-1)), \\
  j &\colon V\cong H^1(\CP^2, E(-2)) \to W\cong H^1(l, E(-2)),
\end{split}\end{equation*}
which are induced from the short exact sequences of sheaves
\begin{gather*}
    0 @>>> \shfO(-1) @>z_0>> \shfO @>>> \shfO_l @>>>0, \\
    0 @>>> \shfO(-3) @>z_0>> \shfO(-2) @>>> \shfO_l(-2) @>>>0.
\end{gather*}

Then we have the following \cite{Vect}:
\begin{Theorem}
  Let $\M(r,0,n)$ be the moduli space of torsion free sheaves $E$ over
  $\CP^2$ with a trivialization $E|_l\cong W\otimes\shfO_l$ which has
  $c_1(E) = 0$, $c_2(E) = n$.
  The above correspondence defines a bijection between $\M(r,0,n)$ and
  the quotient of the set of $(B_1, B_2, i, j)\in\bM$ satisfying
  \begin{enumerate}
  \item \rom(complex ADHM equation\rom) $\mu_\C(B_1, B_2, i, j) =
    [B_1, B_2] + ij = 0$
  \item $\sigma$ is injective and $\tau$ is surjective,
  \end{enumerate}
  by the natural action of $\GL(V)$. Here the Chern classes are given
  by $c_1(E) = 0$, $c_2(E) = \dim V$.
\label{thm:monad}\end{Theorem}

We can replace the second condition by what is similar to one
introduced in \cite[3.5]{Na-quiver},\cite[3.8]{Na-alg}.
(In fact, the condition used there is the `transpose' of the following.)
The proof is easy and omitted.
\begin{Lemma}
  The injectivity of $\sigma$ holds always and the surjectivity of
  $\tau$ is equivalent to saying that there exists no proper subspace
  $S$ of $V$ such that $B_k(S)\subset S$ \rom($k=1,2$\rom) and
  $i(W)\subset S$.
\label{lem:stable}\end{Lemma}

Now we put hermitian metrics on $V$ and $W$. Then $\bM$ is regarded as
a K\"ahler manifold and we can consider the real moment map with
respect to the natural action by $\U(V)$:
\begin{equation*}
  \mu_\R(B_1, B_2, i, j) \defeq [B_1, B_1^\dagger] + [B_2, B_2^\dagger]
  + i i^\dagger - j^\dagger j,
\end{equation*}
where $(\cdot)^\dagger$ is the hermitian adjoint.

Thanks to the above
lemma, the following can be proved exactly as in
\cite[3.5]{Na-quiver}:
\begin{Proposition}
  Fix a negative real number $\zeta_\R$. Then the
  condition~$(2)$ in \thmref{thm:monad} holds if and only if there
  exists $g\in\GL(V)$ such that
  \begin{equation*}
    \mu_\R(g B_1 g^{-1}, g B_2 g^{-1}, g i, j g^{-1}) = - \zeta_\R.
  \end{equation*}
Moreover such a $g$ is unique up to $\U(V)$.
\end{Proposition}

Hence the moduli space can be identified with
\begin{equation*}
  \mu_\C^{-1}(0)\cap\mu_\R^{-1}(-\zeta_\R) / \U(V).
\end{equation*}
The map $\mu = (\mu_\R, \mu_\C)$ is a hyper-K\"ahler moment map in the
sense of \cite{HKLR},
where the vector space $\bM$ has a structure of a $\Bbb H$-module by
\begin{equation*}
  J (B_1, B_2, i, j) \defeq
  (-B_2^\dagger, B_1^\dagger, -j^\dagger, i^\dagger).
\end{equation*}
Thus the moduli space $\M(r,0,n)$ is a hyper-K\"ahler quotient.

Since the action of $\U(V)$ on
$\mu_\C^{-1}(0)\cap\mu_\R^{-1}(-\zeta_\R)$ is free, we have
(cf.~\cite[1.3]{Na-resol})
\begin{Theorem}
  The moduli space $\M(r,0,n)$ is a smooth hyper-K\"ahler manifold of
  real dimension $4nr$. Moreover the metric is complete.
\end{Theorem}

We would like to remark the relation with the original ADHM
description. This was already explained in \cite{Na-resol} in more
detail. So our description is sketchy.
The framed moduli space of
$\SU(r)$-instanton with $c_2 = n$ is described as a hyper-K\"ahler
quotient
\begin{equation}
  \{ (B_1, B_2, i, j) \in \mu_\C^{-1}(0)\cap\mu_\R^{-1}(0) \mid
  \text{stabilizer in $\U(V)$ is trivial}\, \} / \U(V)\, .
\label{eq:instanton}\end{equation}
This is also a smooth hyper-K\"ahler manifold, but the metric is {\it
  not\/} complete. The metric completion
$\mu_\C^{-1}(0)\cap\mu_\R^{-1}(0)/\U(V)$ is isomorphic to the framed
moduli space of ideal instantons, or Uhlenbeck's (partial)
compactification. This space has singularities, and there exists a
natural morphism from our space $\M(r,0,n)$, which is a resolution.
These two moduli spaces, of torsion free sheaves and of ideal
instantons, naturally appear for more general projective surfaces. The
construction of the morphism for general cases was done by J.~Li
\cite{Li}.

Suppose $r = 1$. The double dual $E^{\vee\vee}$ of a rank
$1$ torsion-free sheaf $E$ is locally free and has
$c_1(E^{\vee\vee}) = 0$. Hence $E^{\vee\vee} = \shfO$ and $E$ is an
ideal $\cal J$ of colength finite. This shows that $\M(1,0,n)$ is
isomorphic to the Hilbert scheme $\Hilbn{(\C^2)}$.

\begin{Remark}
(1)\
Note that the $\Hilbn{K3}$ was the first example of higher dimensional
compact hyper-K\"ahler manifold given by Fujiki and Beauville~\cite{Be}.
It seems natural to conjecture that $\Hilbn{X}$ has a hyper-K\"ahler
structure if $X$ has a hyper-K\"ahler structure. For example,
when $X$ is an ALE space, this is true.
When $\Hilbn{X}$ is projective, the existence of a holomorphic
symplectic structure is enough thanks to the solution of the Calabi
conjecture by Yau~\cite{Yau}.
Although there are some extensions of Calabi conjectures to noncompact
manifolds (e.g., \cite{BK,TY}), our case $\Hilbn{(\C^2)}$ is uncovered
probably. This is because our manifolds do not have quadratic curvature
decay, and the cones of the ends are not smooth manifolds.

(2)\
In \cite{KN}, the ADHM description for instantons on ALE spaces was
given. This is similar to \eqref{eq:instanton}, but the vector space
$V$ now becomes a representation of quivers of affine type.
The modification ``$\mu_\R^{-1}(0) \to \mu_\R^{-1}(-\zeta_\R)$'' used
in this appendix can be also applied to the case of ALE spaces.
It gives the ADHM description of torsion free sheaves, instead of
instantons.
In fact, the corresponding moduli spaces were already studied in
\cite{Na-quiver}.
More precisely, the parameter $\zeta$ used in \cite{KN} was tracefree,
while one used in \cite{Na-quiver} was not.
Thus the action of affine Lie algebra, constructed in \cite{Na-quiver},
is not on the homology of moduli spaces of instantons, but of torsion
free sheaves.
In \cite{Na-gauge}, we persisted in {\it genuine\/} vector bundles and
we give only representations of {\it finite\/} dimension
Lie algebras.

(3)\
There is yet another modification of \eqref{eq:instanton}.
In stead of changing the value of the real moment map $\mu_\R$, we can
change the value of the complex moment map to get
$\mu_\R^{-1}(0)\cap\mu_\C^{-1}(-\zeta_\C)/\U(V)$.
It is a deformation of $\mu_\R^{-1}(-\zeta_\R)\cap\mu_\C^{-1}(0)/\U(V)$.
Changing the identification $\R^4 \cong \C^2$, we can identify
this space also with the moduli space of torsion free sheaves.
\label{sheaf}\end{Remark}
\subsection{Topology of $\Hilbn{(\C^2)}$}
Our next task is to compute the homology group of the Hilbert scheme
$\Hilbn{(\C^2)}$.

We consider $\Hilbn{(\C^2)}$ as a hyper-K\"ahler quotient as in the
previous section.
As in \cite{ES,Na-homology}, the idea is to use the torus action given
by
\begin{equation*}
  (B_1, B_2, i, j) \mapsto (t_1 B_1, t_2 B_2, t_1 i, t_2 j)
  \qquad (t_1, t_2)\in S^1\times S^1.
\end{equation*}
This action commutes with the $\U(V)$-action and induces an action on
$\M(1,0,n)$.
If $\U(V).(B_1, B_2, i, j)\in\M(1,0,n)$ is a fixed point of the torus
action, there exists a homomorphism $\lambda\colon\C^\ast\to\GL(V)$
such that
\begin{equation*}\begin{split}
  & (t_1 B_1, t_2 B_2, t_1 i, t_2 j) \\
  =\; &(\lambda(t_1,t_2) B_1 \lambda(t_1,t_2)^{-1},
     \lambda(t_1,t_2) B_2 \lambda(t_1,t_2)^{-1}, \lambda(t_1,t_2)i,
     j\lambda(t_1,t_2)^{-1}).
\end{split}\end{equation*}
If $V = \bigoplus V(k,l)$ is the weight space decomposition, the above
equation can be written as
\begin{gather*}
  B_1(V(k,l))\subset V(k+1,l), \quad B_2(V(k,l))\subset V(k,l+1), \\
  i(W)\subset V(1,0), \quad j(V(k,l)) = 0 \;\text{unless $(k,l) = (0,-1)$}.
\end{gather*}

Since $\bigoplus_{k \ge 1, l \ge 0} V(k,l)$ contains $\Ima i$ and
$B_k$-invariant, the condition in \lemref{lem:stable} implies
it is equal to $V$. In particular, $j = 0$.
The situation is visualized as
\begin{equation*}
\begin{CD}
  \C = W @>i>> V(1,0) @>B_1>> V(2,0) @>B_1>> V(3,0) @>B_1>> \cdots\\
         @.   @VVB_2V         @VVB_2V        @VVB_2V        \\
         @.    V(1,1) @>B_1>> V(2,1) @>B_1>> V(3,1) @>B_1>> \cdots\\
         @.   @VVB_2V         @VVB_2V        @VVB_2V        \\
         @.    \vdots @.      \vdots @.      \vdots
\end{CD}
\end{equation*}

The easy argument using the induction on $k$ or $l$ shows
\begin{Lemma}
  \begin{enumerate}
  \item   $\dim V(k,l) = 1$ or $0$,
  \item   $\sum_l \dim V(k,l)$ is monotone non-increasing with respect
    to $k$.
  \end{enumerate}
\end{Lemma}

The fixed point is uniquely determined by $\dim V(k,l)$, and hence by
the partition
\begin{equation*}
  (\sum_l \dim V(1,l), \sum_l \dim V(2,l), \dots)
\end{equation*}
of $n$.

The map
\begin{equation*}
  F(\U(V).(B_1, B_2, i, j)) \defeq
  \|B_1\|^2 + \varepsilon \|B_2\|^2 + \| i\|^2 + \varepsilon \|j\|^2
\end{equation*}
is the moment map for the torus action coupled with a certain element
in the Lie algebra of the torus.
If we take $\varepsilon > 0$ sufficiently small and generic,
the critical point $p$ of $F$ is a fixed point.
Then the tangent space $T_p\Hilbn{(\C^2)}$ is a torus module and
the Morse index of $F$ at $p$ is given by the sum of the dimension of
weight spaces $T(k,l)$ such that
(a) $k < 0$, or (b) $k = 0$ and $l < 0$.
Calculating the dimension, we find
\begin{equation*}
  \text{Index of $F$ at $p$}
  = 2 \sum_{l \ge 1} \sum_k \dim V(k,l)
  = 2n - 2\sum_k \dim V(k,0).
\end{equation*}
If $p$ corresponds to a partition of $n$ as above, this index is equal
to
\begin{equation*}
  2n - 2(\text{number of parts}).
\end{equation*}

Thus we get
\begin{Theorem}
  The homology group $H_*(\Hilbn{(\C^2)})$ has no torsion and vanishes
  in odd degrees.
  The Betti number of $\Hilbn{(\C^2)}$ is given by
  \begin{equation*}
    b_{2i}(\Hilbn{(\C^2)}) = \text{number of partitions of $n$ into
      $n - i$ parts}.
  \end{equation*}
\end{Theorem}

\begin{Corollary}
  If $P_t(\Hilbn{(\C^2)}) = \sum t^i b_i(\Hilbn{(\C^2)})$ denotes the
  Poincar\`e polynomial, its generating function is given by
  \begin{equation*}
    \sum_{n\ge 0} q^n P_t(\Hilbn{(\C^2)})
    = \prod_{m=1}^\infty \frac 1{1 - t^{2m - 2}q^m}.
  \end{equation*}
\end{Corollary}

\end{document}